\documentclass[10pt,twocolumn,a4paper]{article}

\usepackage[T1]{fontenc}
\usepackage[utf8]{inputenc}
\usepackage{mathptmx}            
\usepackage[scaled=0.92]{helvet} 
\usepackage{microtype}

\usepackage[a4paper,top=2.2cm,bottom=2.4cm,left=1.8cm,right=1.8cm,columnsep=0.7cm]{geometry}
\usepackage{float}
\usepackage{amsmath,amssymb}

\usepackage{graphicx}

\usepackage[table]{xcolor}
\usepackage{booktabs}
\usepackage{caption}
\usepackage{subcaption}
\usepackage{titlesec}
\usepackage{fancyhdr}
\usepackage{enumitem}
\usepackage[hidelinks,breaklinks]{hyperref}
\usepackage{balance}

\definecolor{psrblue}{HTML}{1F3864}
\definecolor{accent}{HTML}{2E75B6}
\definecolor{rulegray}{HTML}{BFBFBF}

\titleformat{\section}
  {\sffamily\bfseries\color{psrblue}\large}{\thesection}{0.6em}{\MakeUppercase}
\titleformat{\subsection}
  {\sffamily\bfseries\color{psrblue}\normalsize}{\thesubsection}{0.5em}{}
\titlespacing*{\section}{0pt}{1.4ex plus 1ex minus .2ex}{0.8ex plus .2ex}
\titlespacing*{\subsection}{0pt}{1.1ex plus 1ex minus .2ex}{0.5ex plus .2ex}

\fancypagestyle{main}{%
  \fancyhf{}

  \fancyhead[L]{\footnotesize\sffamily\color{psrblue}HERA-S: Closed-Loop Pumped-Storage Hydropower in Rio de Janeiro}
  \fancyhead[R]{\footnotesize\sffamily\color{psrblue}\thepage}
}
\fancypagestyle{firstpage}{%
  \fancyhf{}

  \fancyfoot[C]{\scriptsize\thepage}
}

\newcommand{\thetitle}{Optimization of Closed-Loop Pumped-Storage Hydropower Siting}

\begin{document}

\twocolumn[{%
\begin{center}
  {\sffamily\bfseries\LARGE\color{psrblue} \thetitle \par}
  \vspace{1.1em}
  {\large
   Luiz Rodolpho Albuquerque\textsuperscript{1}\quad
   Rafael Kelman\quad
   Tarcisio Castro\quad
   Ana Petrungaro\quad
   Tiago Andrade\par}
  \vspace{0.5em}
  {\itshape PSR, Rio de Janeiro, Brazil\par}
  \vspace{0.55em}
\end{center}

\vspace{0.4em}
\begin{center}
\begin{minipage}{0.92\textwidth}
  \color{black}
  {\small
  \noindent\textbf{\textsf{\color{psrblue}ABSTRACT}}\\[0.3em]
  Accelerating variable renewable energy integration introduces substantial operational challenges to modern power grids. Frequent curtailment of surplus renewable generation highlights the pressing need for long-duration energy storage capable of shifting generation to high-demand periods while providing essential ancillary services and grid inertia. Socio-environmental constraints increasingly limit conventional hydropower development; however, closed-loop pumped-storage hydropower (PSH) presents a viable alternative because its off-river reservoirs avoid natural watercourses. This paper introduces HERA-S, a computational modeling framework developed by PSR (supported by EDF, CTG, Brookfield, and Rio Light) to streamline and standardize regional PSH site prospecting. HERA-S automates spatial screening, dam optimization, cost estimation, and socio-environmental impact scoring. Applied to a case study in Rio de Janeiro, Brazil, the framework employs an integer programming model to optimize dam-fill and excavation mass balances against inter-reservoir distance, significantly improving pre-feasibility planning agility.\\[0.6em]

  \textbf{Keywords:} Pumped-storage hydropower; geoprocessing; integer optimization; digital elevation model; water-resources planning; variable renewable energy; energy sustainability.}
\end{minipage}
\end{center}
\vspace{1.3em}
}]

\thispagestyle{firstpage}

\section{Introduction}

The global energy transition introduces complex balancing challenges for regional power transmission grids. As variable renewable energy (VRE) resources---principally wind and solar photovoltaics---expand, supply-demand volatility increases. Modern power systems require flexible, large-scale storage resources to absorb off-peak surpluses and dispatch energy during shortages. Pumped-storage hydropower (PSH) remains the most commercially mature and cost-effective technology for utility-scale energy storage.

PSH operates by pumping water from a lower reservoir to an upper elevation during periods of low net demand (power consumption minus VRE generation), storing gravitational potential energy. During peak demand or sudden generation drops (e.g., solar ramps at dusk), the water generates electricity through hydro-turbines. Compared to competing technologies, PSH offers long operational lifespans, high round-trip efficiency, fast dynamic response, and dynamic voltage/frequency support. Furthermore, its incremental energy storage cost is lower than lithium-ion battery energy storage systems (BESS) for long durations \cite{ref1}, even accounting for recent battery cost declines.

The main economic advantage of long-duration PSH stems from the structural decoupling of power output and energy capacity. Power capability is constrained by fixed electromechanical infrastructure (powerhouses, hydro-turbines, and conveyance tunnels), whereas energy storage capacity scales primarily with reservoir volume. Expanding storage capacity relies mainly on raising dam heights or rebalancing site earthworks rather than expanding the hydraulic machinery. Consequently, the marginal cost per added kWh of storage capacity declines rapidly as duration increases.

Conversely, BESS scales modularly and linearly. Because electrochemical cells comprise both the storage medium and energy conversion mechanism, doubling capacity requires nearly double the physical battery modules. Adding BESS duration exhibits a constant marginal cost per kWh. Furthermore, while PSH inherently supplies physical mechanical inertia to support system frequency, BESS requires advanced grid-forming inverters, adding capital cost.

Thus, PSH marginal costs decrease significantly with longer discharge durations ($>8$\,h) wherever favorable topography exists, while battery scalability remains constrained by cell manufacturing economics. 

In Brazil, interest in PSH is expanding to complement the country's extensive conventional hydro fleet amid massive wind growth in the Northeast and widespread solar buildout nationwide. Although Brazilian topography offers promising sites (e.g., Serra da Mantiqueira, Serra do Mar, and the Southern Plateau), systematic PSH prospecting remains limited. Existing workflows lack integrated methodologies that combine geospatial terrain evaluation, automated civil engineering costing, socio-environmental screening, and regional portfolio optimization. The HERA framework directly addresses this methodological gap.

\subsection{Background}
Conventional hydropower planning in Brazil historically followed the \emph{Manual for the Hydropower Inventory of River Basins} \cite{ref2}, evaluating basins via field reconnaissance, local alternative comparisons, and project-by-project environmental reviews. While effective historically, this approach presents two major limitations in modern planning:
\begin{enumerate}[leftmargin=*]
  \item It evaluates restricted pre-selected dam axes without systematically scanning all topographically viable basin sites, risking suboptimal selections.
  \item It evaluates project impacts in isolation, missing cumulative ecological impacts such as river fragmentation or basin connectivity loss.
\end{enumerate}

HERA overcomes these limitations by conducting systematic spatial screenings across entire basins while integrating cumulative socio-environmental constraints directly into economic head-optimization models.

For closed-loop PSH, these considerations differ. The standard layout (two reservoirs separated by elevation, linked via penstocks and a reversible powerhouse) requires high-resolution spatial analysis of Digital Elevation Models (DEMs) to identify site pairs satisfying gross head, storable volume, and distance constraints while minimizing civil works (dams, earthworks, and penstock runs).

\subsection{HERA-S Architecture}
HERA was originally created to evaluate conventional basin-scale hydropower potential by integrating: (i) geospatial analysis with socio-environmental layers; (ii) automated civil engineering parametric costing; and (iii) mathematical portfolio optimization under environmental constraints \cite{ref3}.

The extension presented here, \textbf{HERA-S} (\emph{storage}), adapts this architecture for PSH. Its updated engineering module evaluates specialized PSH geometries and solves a mixed-integer optimization problem to identify the minimum-cost artificial reservoir configuration for a targeted head, capacity, and stored energy.

Unlike conventional cascading hydro projects, closed-loop PSH plants operate independently from natural river flows. Each candidate site can be sized independently across various installed capacities (MW) and storage durations (MWh).

HERA-S models multiple cycle configurations. Semi-open systems use an existing river, reservoir, or lake for one basin, reducing capital expenditure (CAPEX) but introducing external hydrological constraints. Closed-loop systems construct two new off-river reservoirs without connecting to natural waterways. This configuration minimizes ecological interference, simplifies environmental permitting, and maximizes operational autonomy.

\begin{figure*}[htbp]
  \centering
  \begin{subfigure}{0.485\textwidth}\includegraphics[width=\linewidth]{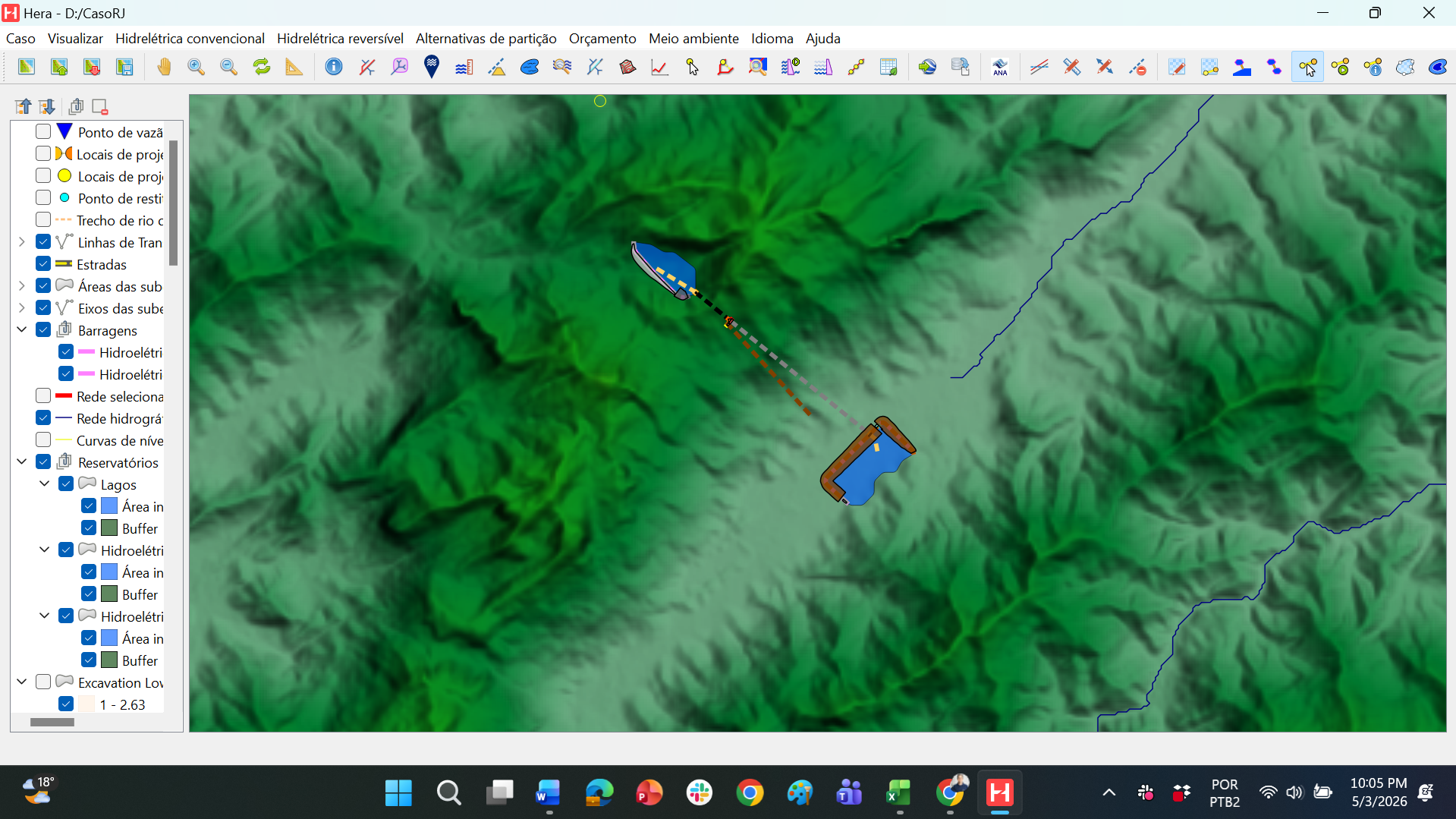}\end{subfigure}\hfill
  \begin{subfigure}{0.485\textwidth}\includegraphics[width=\linewidth]{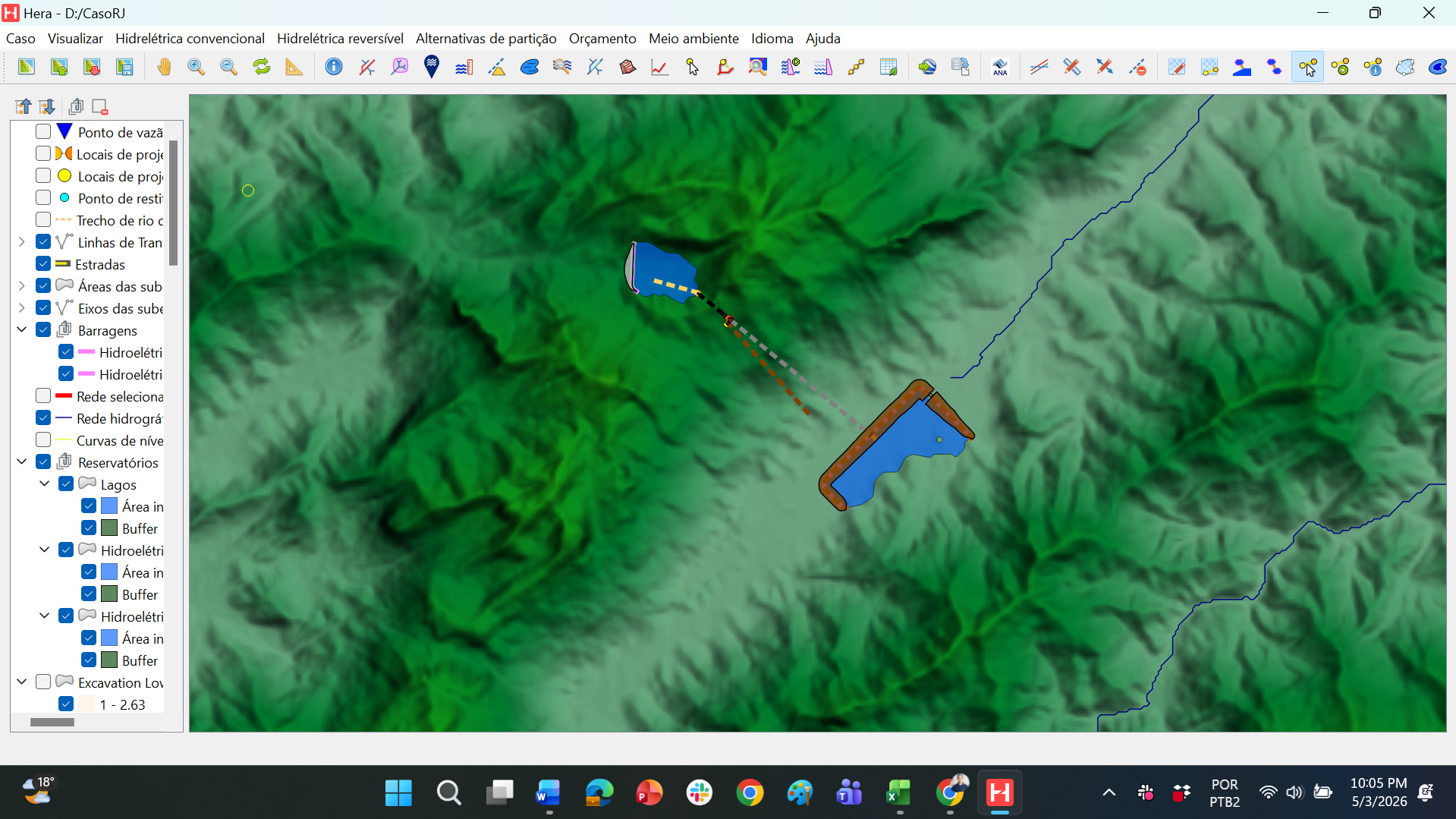}\end{subfigure}

  \vspace{0.5em}
  \begin{subfigure}{0.485\textwidth}\includegraphics[width=\linewidth]{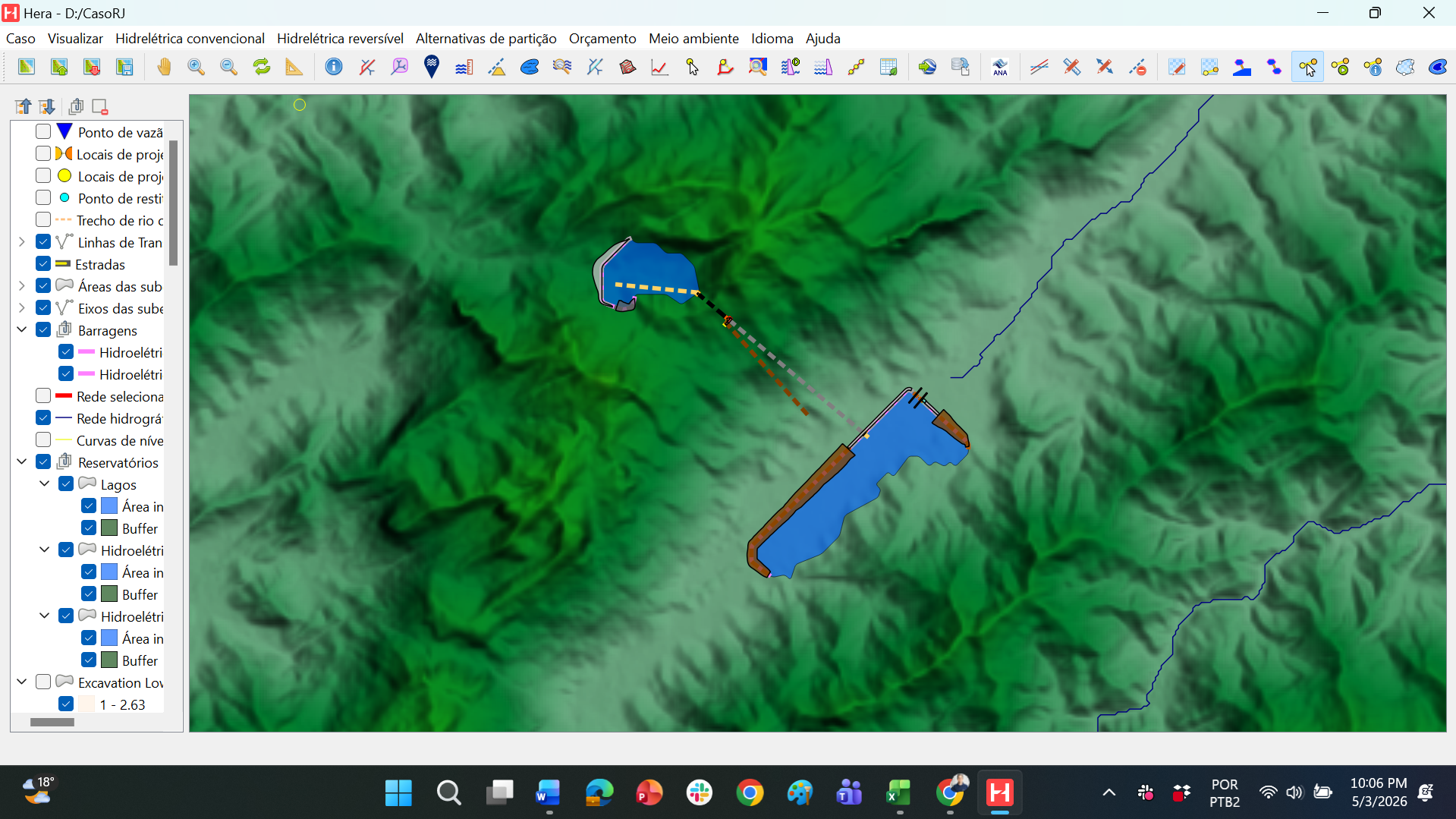}\end{subfigure}\hfill
  \begin{subfigure}{0.485\textwidth}\includegraphics[width=\linewidth]{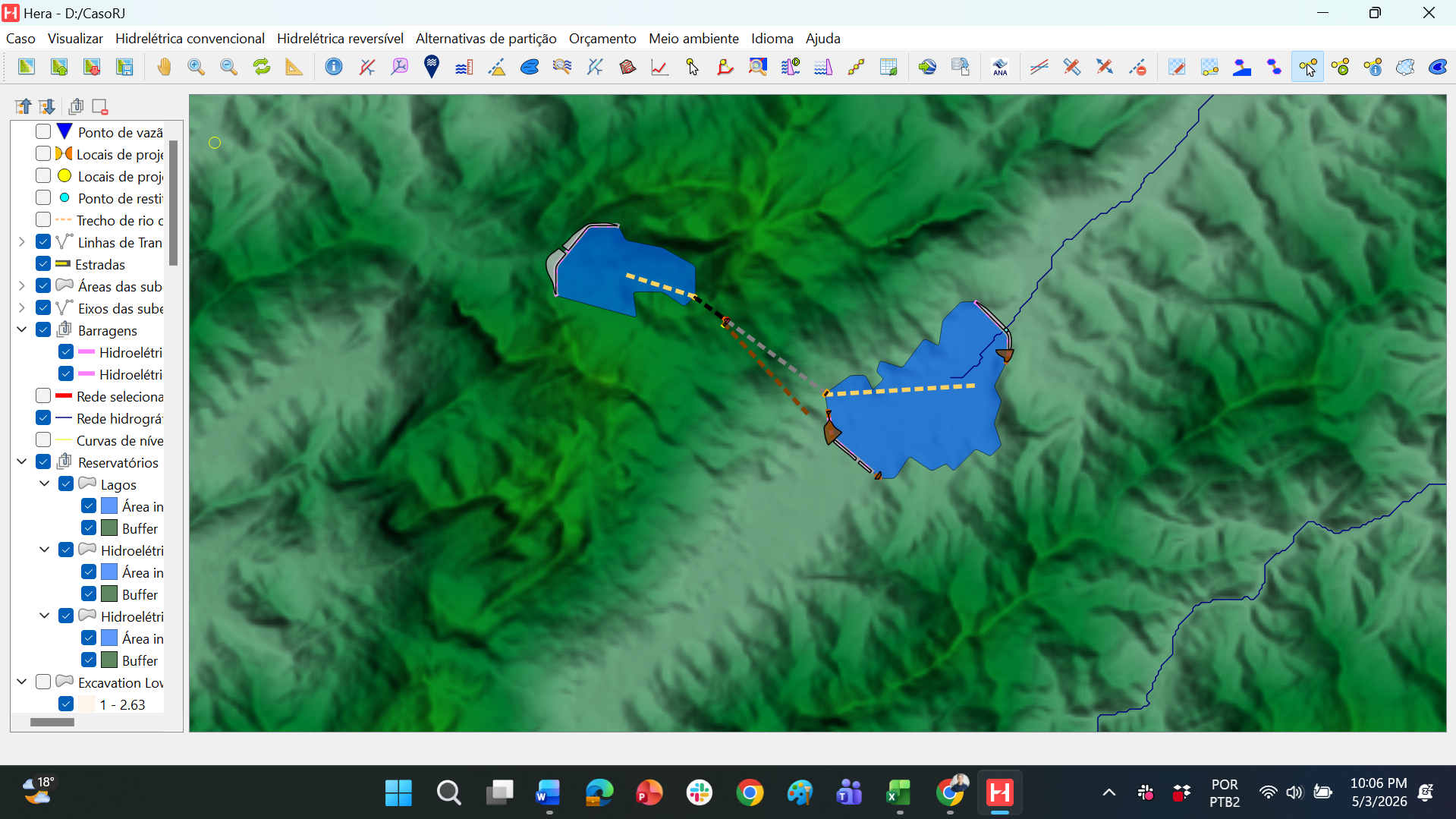}\end{subfigure}

  \vspace{0.5em}
  \begin{subfigure}{0.485\textwidth}\includegraphics[width=\linewidth]{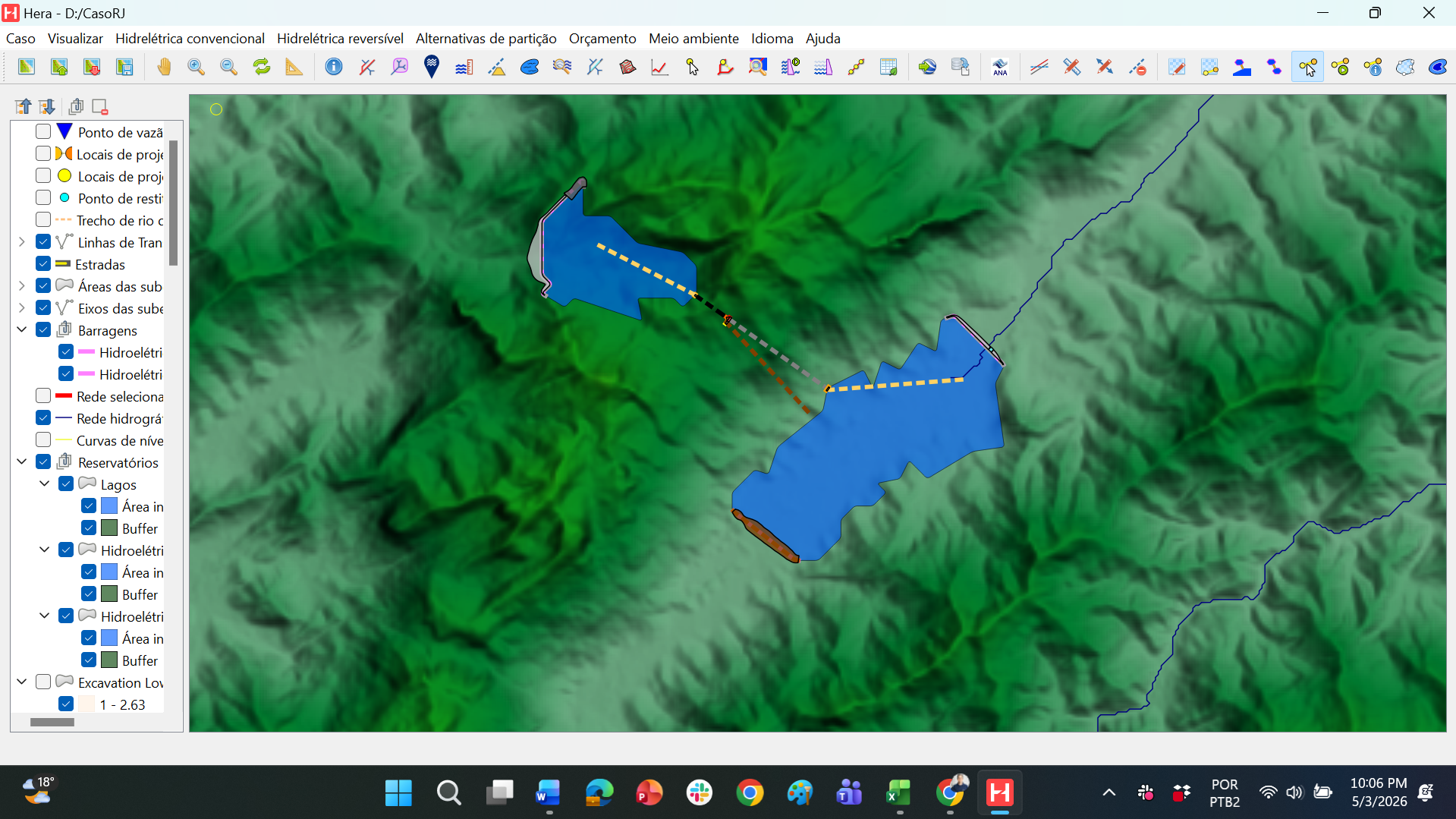}\end{subfigure}\hfill
  \begin{subfigure}{0.485\textwidth}\includegraphics[width=\linewidth]{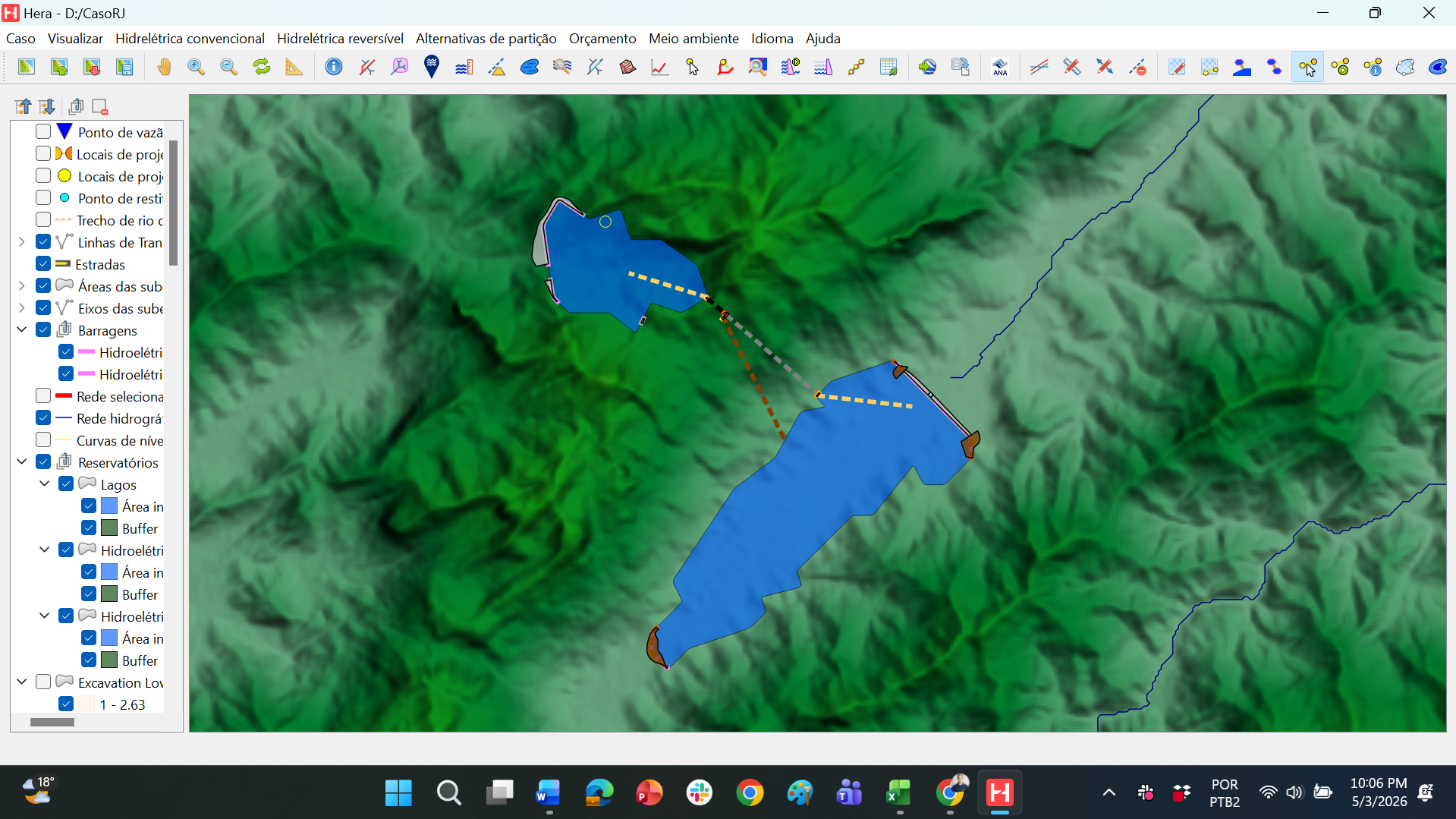}\end{subfigure}
  \caption{Representative 1\,GW PSH plant configurations generated by HERA-S for 6\,h, 12\,h, 24\,h, 48\,h, 72\,h, and 120\,h storage durations.}
  \label{fig:fig1}
\end{figure*}

Figure~\ref{fig:fig1} illustrates HERA-S outputs for six different storage duration configurations (6\,h to 120\,h) at a single 1\,GW candidate site.

\begin{figure*}[htbp]
  \centering
  \includegraphics[width=0.9\textwidth]{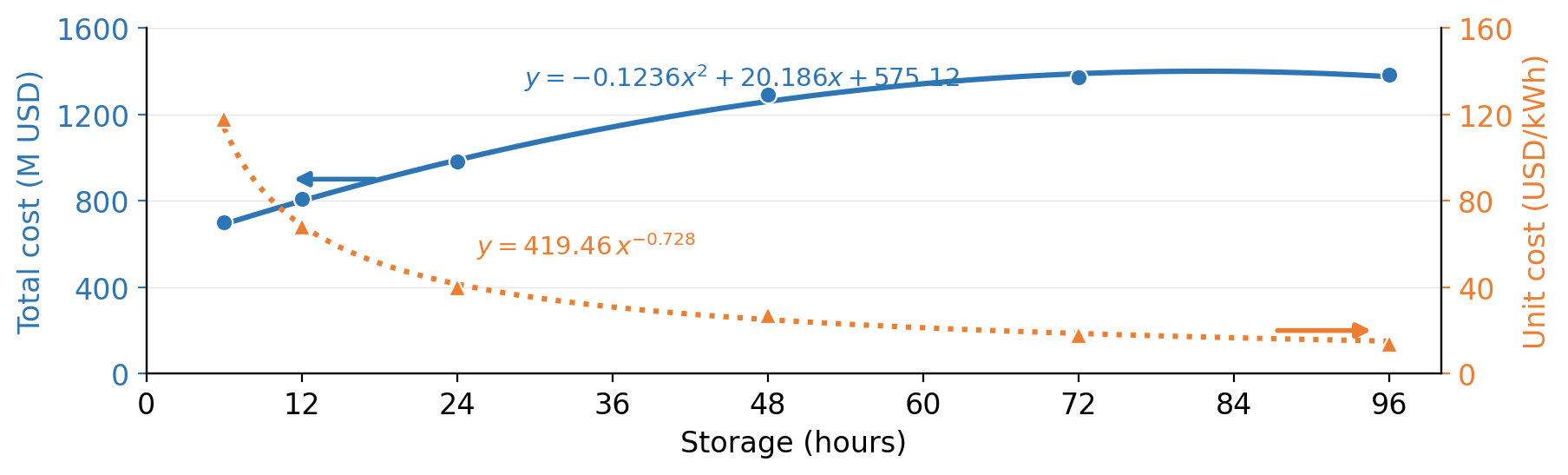}
  \caption{Total capital investment (USD, primary axis) and unit storage cost (USD/kWh, secondary axis) across discharge durations for the candidate site.}
  \label{fig:fig2}
\end{figure*}

Figure~\ref{fig:fig2} plots the calculated unit storage costs, highlighting PSH economies of scale: unit costs fall below 20\,USD/kWh for multi-day storage configurations ($>72$\,h).

By combining geospatial screening with automated cost estimation, HERA-S generates reliable project inventories across vast regions. These candidate project portfolios can then feed directly into long-term expansion planning tools such as OptGen \cite{ref8}. OptGen evaluates PSH candidates based on broader system benefits---such as peak shaving, transmission deferral, renewable integration, load shifting, and provision of grid stability services (inertia, spinning reserves, and black-start capability).

The remainder of this paper is structured as follows: Section~\ref{sec:method} presents the screening workflow; Section~\ref{sec:opt} outlines the reservoir optimization formulation; Section~\ref{sec:case} details a practical case study in Rio de Janeiro State; and Section~\ref{sec:concl} summarizes key findings.

\section{PSH Screening Pipeline}
\label{sec:method}

\begin{figure*}[htbp]
  \centering
  \includegraphics[width=0.9\textwidth]{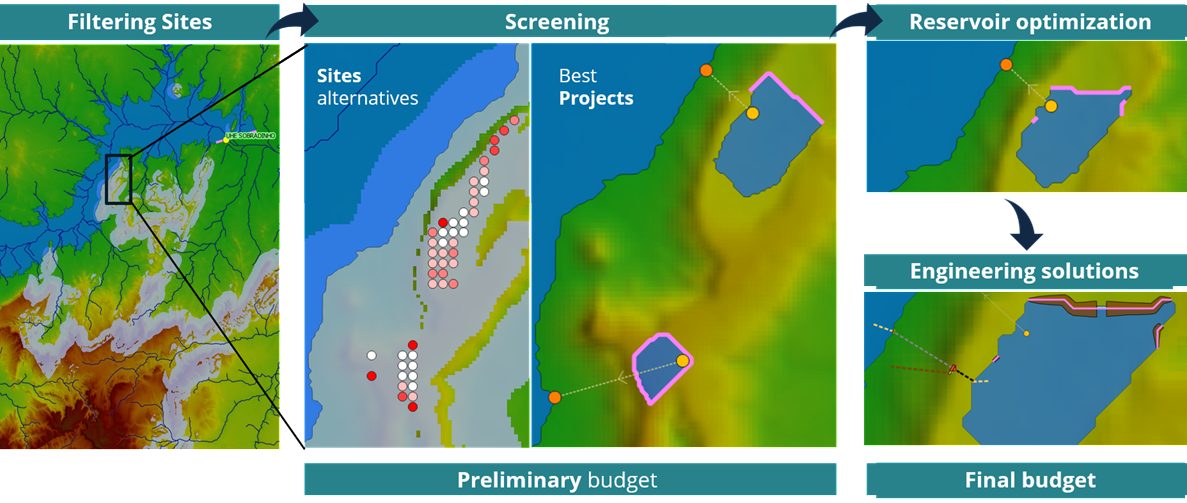}
  \caption{Methodological workflow of HERA-S for PSH site evaluation and optimization.}
  \label{fig:fig3}
\end{figure*}

HERA-S retains the three-stage architecture of original HERA (GIS, Engineering, and Optimization), updated specifically for PSH evaluation (Figure~\ref{fig:fig3}).

The process ingests spatial elevation data alongside thematic constraint layers, applies exclusionary spatial filters, conducts an intensive search for viable site pairs, performs mathematical reservoir optimization, and builds pre-feasibility engineering designs.

\subsection{Filter Application}
The initial stage applies explicit spatial criteria to quickly discard unsuitable land areas. GIS layers define exclusion zones based on user criteria, which constrain the subsequent candidate search space.

Optionally, layers representing protected environmental areas or indigenous lands are overlaid to flag regulatory risks. Excluded areas can be adjusted or re-evaluated dynamically as higher-resolution local data becomes available.

\subsection{Intensive Search}
The intensive search identifies topographically favorable layout pairs within the filtered spatial search space.

The search relies on \emph{geomorphons} \cite{ref7}, an algorithm for fast terrain landform classification using DEMs. User-provided shapefiles can introduce localized geologic, environmental, or land-use constraints. HERA-S then calculates preliminary cost estimates across candidate site combinations, estimating expenditure splits across civil works, equipment, grid interconnections, and environmental mitigation.

\subsection{Engineering Solutions}
Top candidates from the intensive search move to automated structural design and budget refining. HERA-S models civil quantities, electromechanical specifications, tunnel layouts, and site logistics.

The system evaluates reversible Francis or Pelton units for powerhouses, designs underground access tunnels, calculates penstock/conveyance dimensions, and evaluates earthfill, rockfill, or concrete dam options. This generates standardized pre-feasibility budgets across all shortlisted sites.

\section{Reservoir Optimization Model}
\label{sec:opt}

Following candidate identification, HERA-S optimizes artificial reservoir geometries to minimize civil construction costs.

\subsection{Problem Formulation}
Reservoir geometry selection is modeled as a combinatorial integer program over a raster grid of DEM cells. Given a gross head $H$\,(m), installed capacity $P$\,(MW), and operational storage hours $T$\,(h) yielding energy target $P \cdot T$\,(MWh), the objective is to find a cell layout that minimizes total civil expenditure while meeting volume and spatial connectivity requirements.

Following \cite{ref1}, binary variables describe each DEM grid cell: $x_{ij}\in\{0,1\}$ indicates if cell $(i,j)$ forms part of the reservoir boundary (dam embankment or cut wall), and $y_{ij}\in\{0,1\}$ indicates if the cell is within the reservoir interior. Variable $z_{ij}=x_{ij}+y_{ij}$ denotes any cell utilized by the reservoir footprint.

The minimum volume requirement is expressed as:
\begin{equation}
  y_{ij}\,(H-h_{ij})\,\alpha \;\ge\; V_{\min},
\end{equation}
where $h_{ij}$ is the ground elevation of interior cell $(i,j)$, $H$ is the upper reservoir design water level, $\alpha$ is cell area, and $V_{\min}$ is the required active storage volume.

Geometric continuity---preventing fragmented, disconnected water bodies---is enforced via two constraint sets. The first set uses directional separating planes (horizontal, vertical, and diagonal) to block solutions where interior cells bypass perimeter boundaries. If disconnected pockets persist, a second set based on a Traveling Salesman Problem formulation \cite{ref5} eliminates perimeter subtours.

The objective function minimizes total costs for excavation, embankments, slope protection, and high/low-pressure water conveyance conduits. Embankment costs depend on the trapezoidal volume needed to raise boundary cells to elevation $H$. Soil and rock excavation costs are calculated based on overburden depth profiles. Conveyance costs are determined by selecting an optimal boundary connection point that minimizes path length to the lower water body.

\subsection{Solution Strategies}
To solve large grid instances efficiently, HERA-S applies two computational acceleration strategies:
\begin{enumerate}[leftmargin=*]
  \item \textbf{Separating Planes:} Progressively cuts unfeasible topological space without adding decision variables.
  \item \textbf{Progressive Zoom-In Heuristic:} Solves the model first on a coarse, aggregated DEM grid, locates the optimal boundary, clips the DEM around this candidate location, and increases resolution iteratively back to the native DEM scale.
\end{enumerate}

Combining these methods yields optimal or near-optimal solutions within reasonable runtimes for grids up to $266\times266$ cells \cite{ref1}. For larger regional inputs, geospatial pre-screening filters identify promising terrain windows (e.g., natural bowls with high head-to-distance ratios) prior to running the optimization model.

\subsection{Integration with Long-Term Planning}
For each head ($H$) and duration ($\Delta t$) pairing, HERA-S optimizes reservoir dimensions and updates capital cost estimates. These refined outputs provide input data for long-term expansion models like OptGen \cite{ref8}. OptGen determines optimal investment decisions---specifying project locations, commissioning timing, and required storage capacities---to minimize overall power system expansion and operational costs.

\section{Case Study: Rio de Janeiro State}
\label{sec:case}

\subsection{Dataset}
\begin{figure*}[htbp]
  \centering
  \newcommand{\panelw}{0.30\textwidth}
  \begin{subfigure}{\panelw}\includegraphics[width=\linewidth]{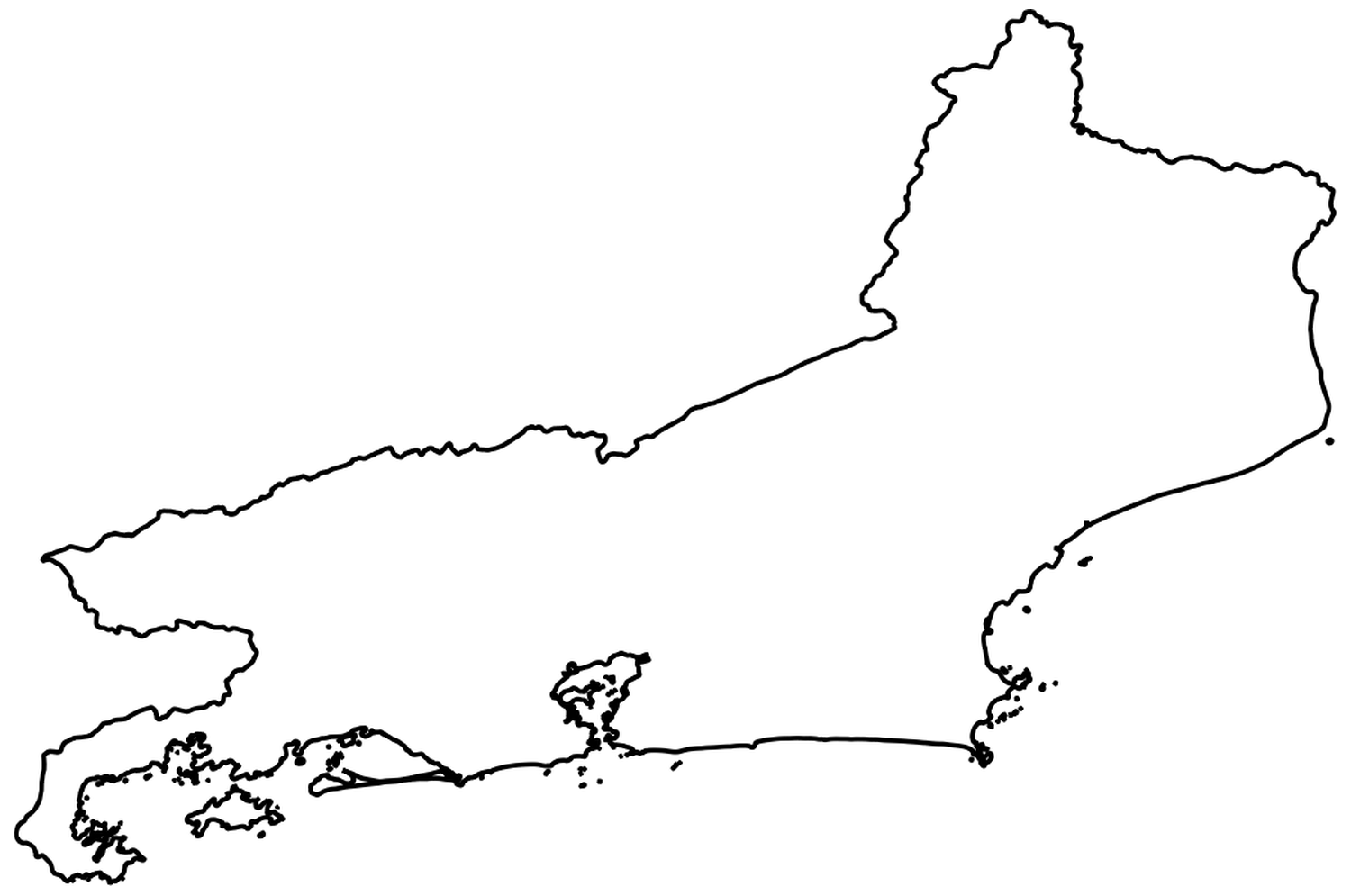}
    \caption*{\footnotesize State boundary}\end{subfigure}\hfill
  \begin{subfigure}{\panelw}\includegraphics[width=\linewidth]{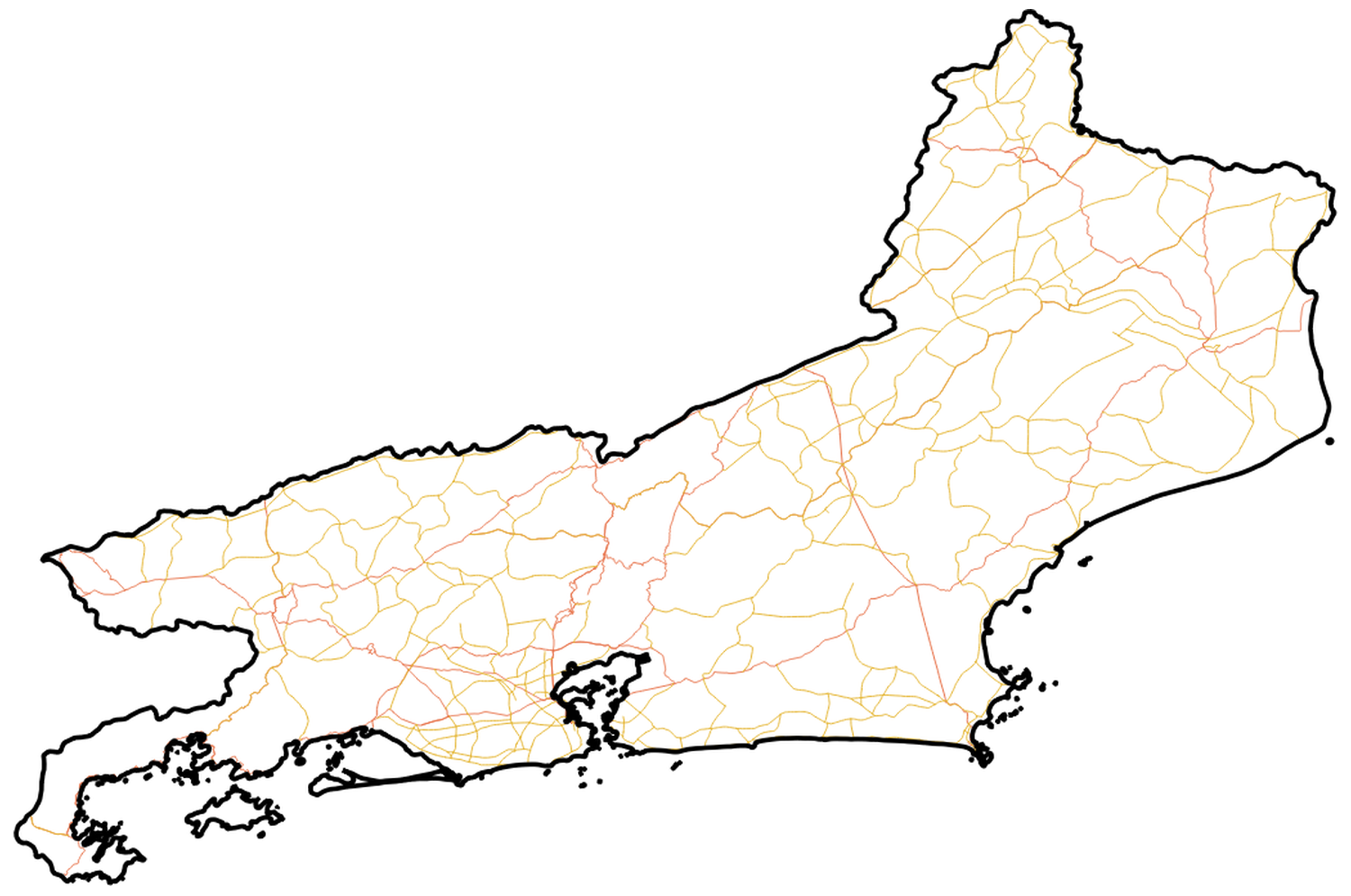}
    \caption*{\footnotesize Transport network}\end{subfigure}\hfill
  \begin{subfigure}{\panelw}\includegraphics[width=\linewidth]{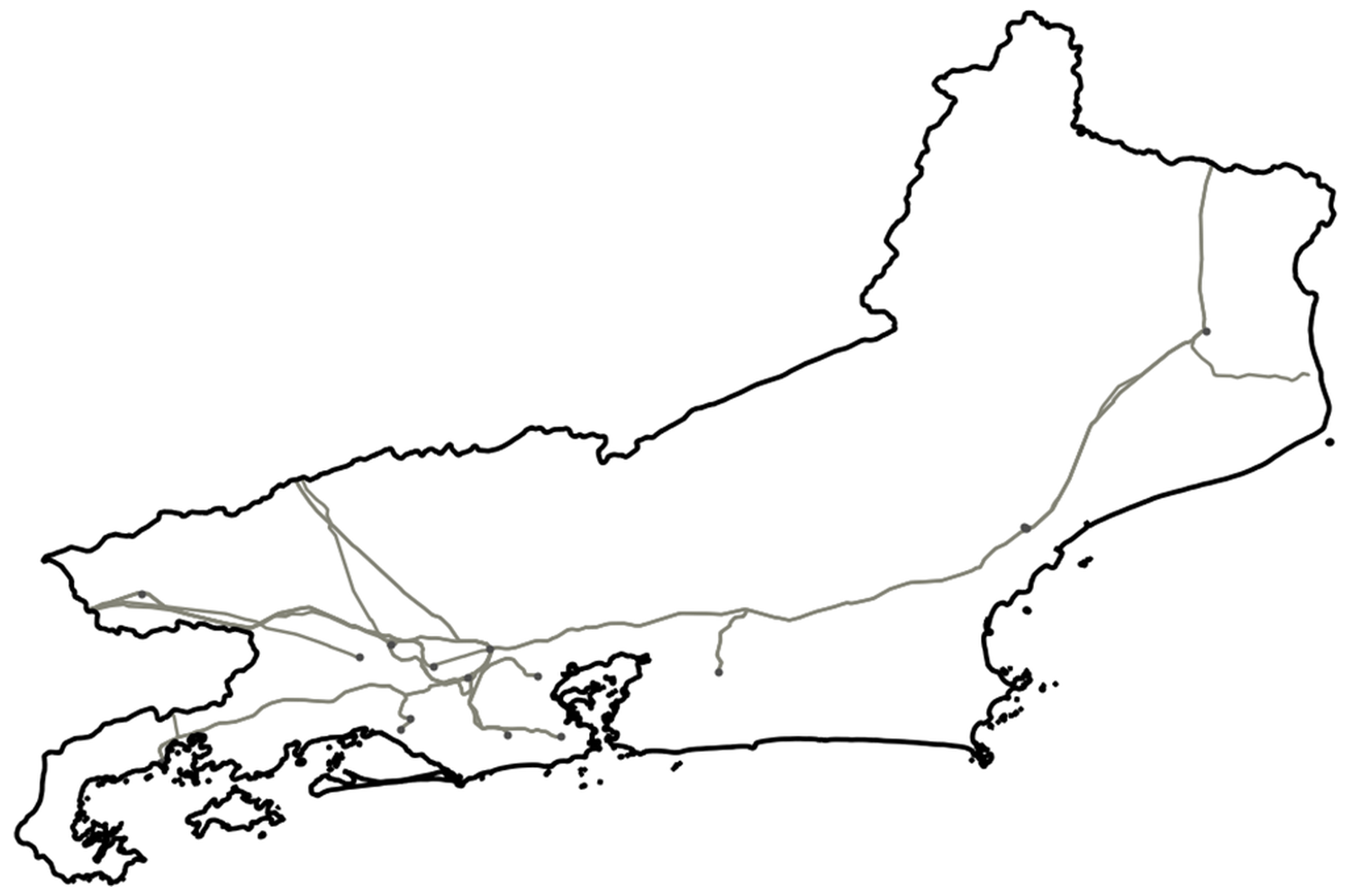}
    \caption*{\footnotesize Grid infrastructure}\end{subfigure}

  \vspace{0.5em}
  \begin{subfigure}{\panelw}\includegraphics[width=\linewidth]{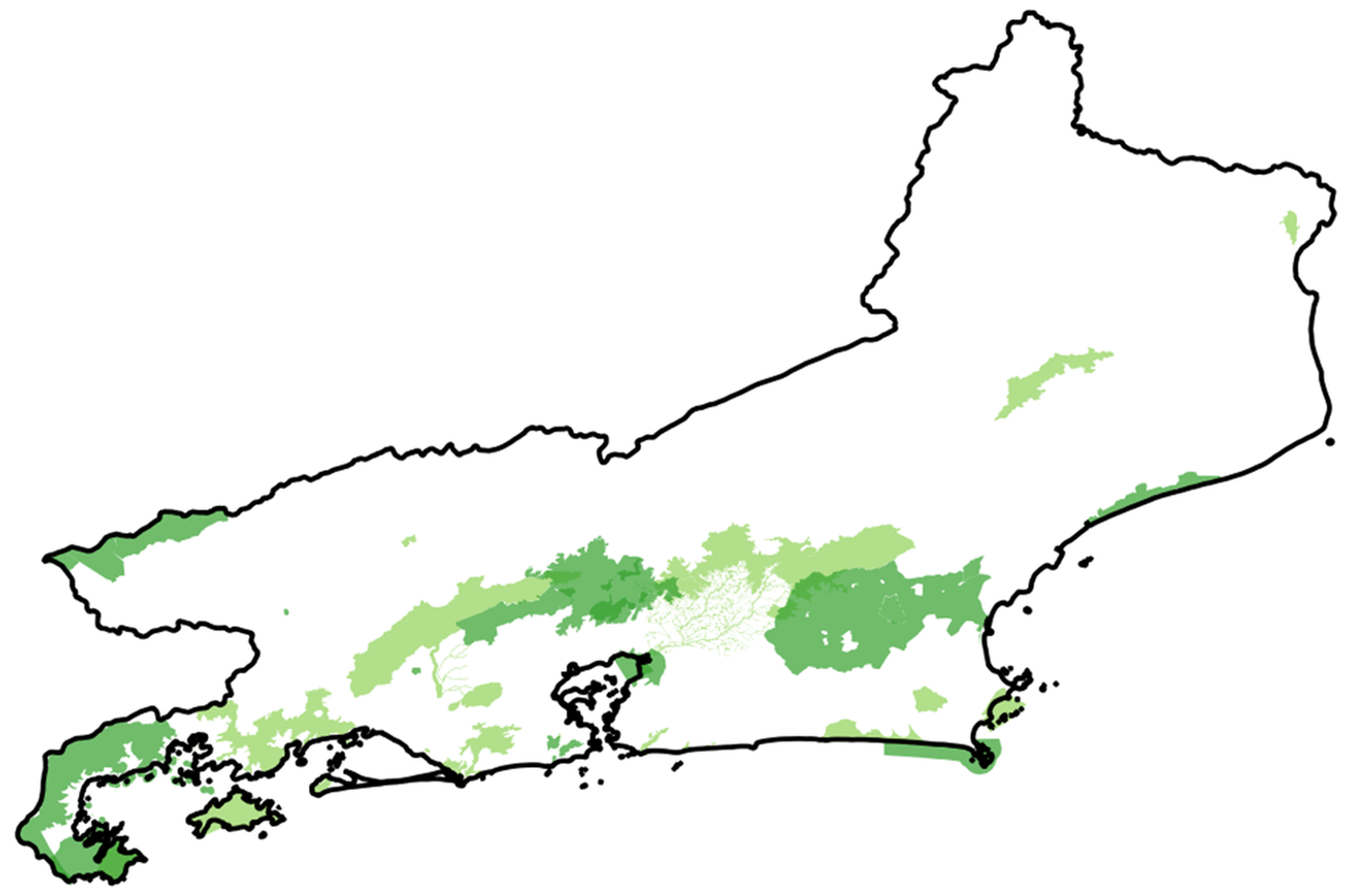}
    \caption*{\footnotesize Conservation units}\end{subfigure}\hfill
  \begin{subfigure}{\panelw}\includegraphics[width=\linewidth]{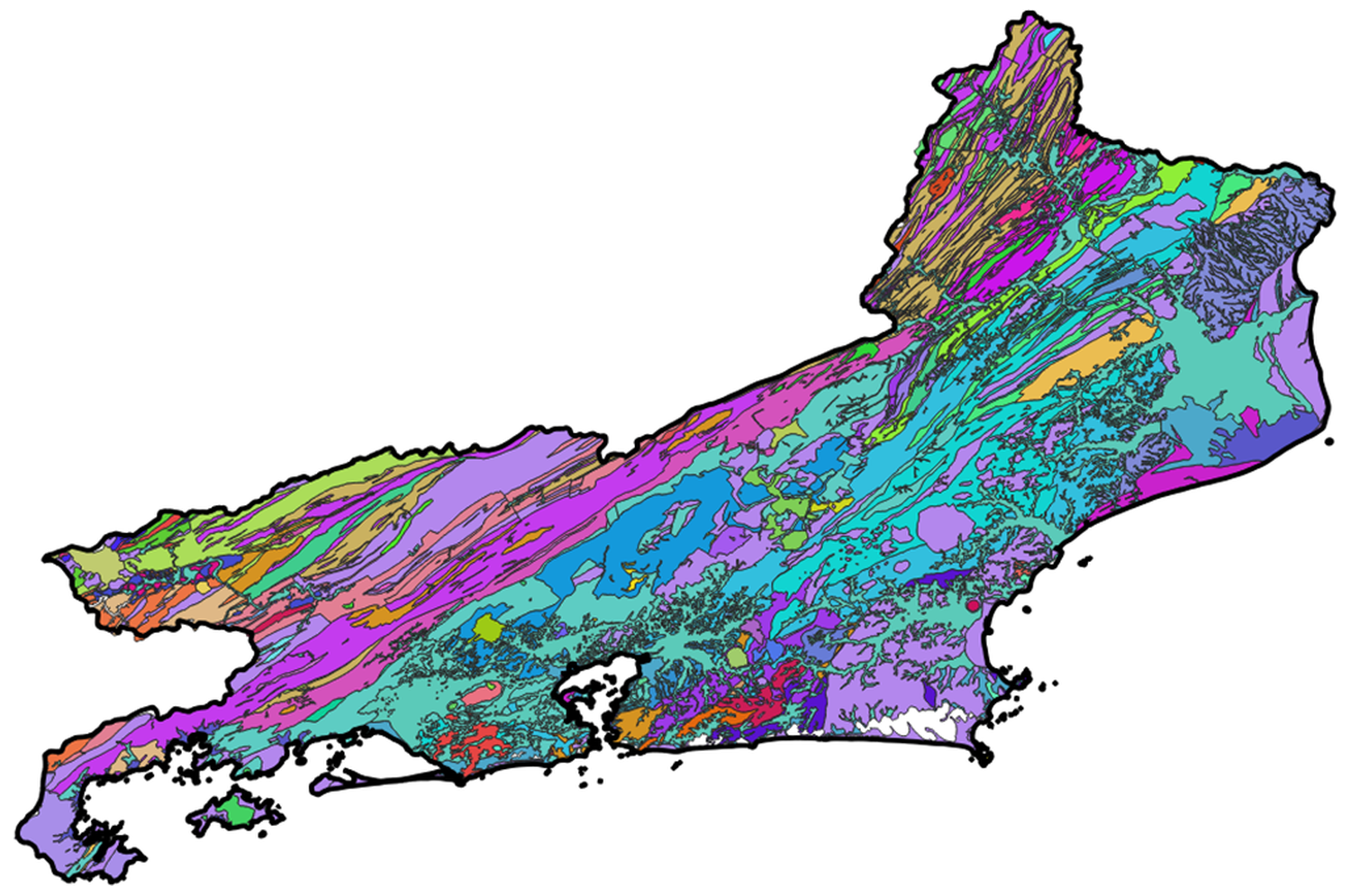}
    \caption*{\footnotesize Regional geology}\end{subfigure}\hfill
  \begin{subfigure}{\panelw}\includegraphics[width=\linewidth]{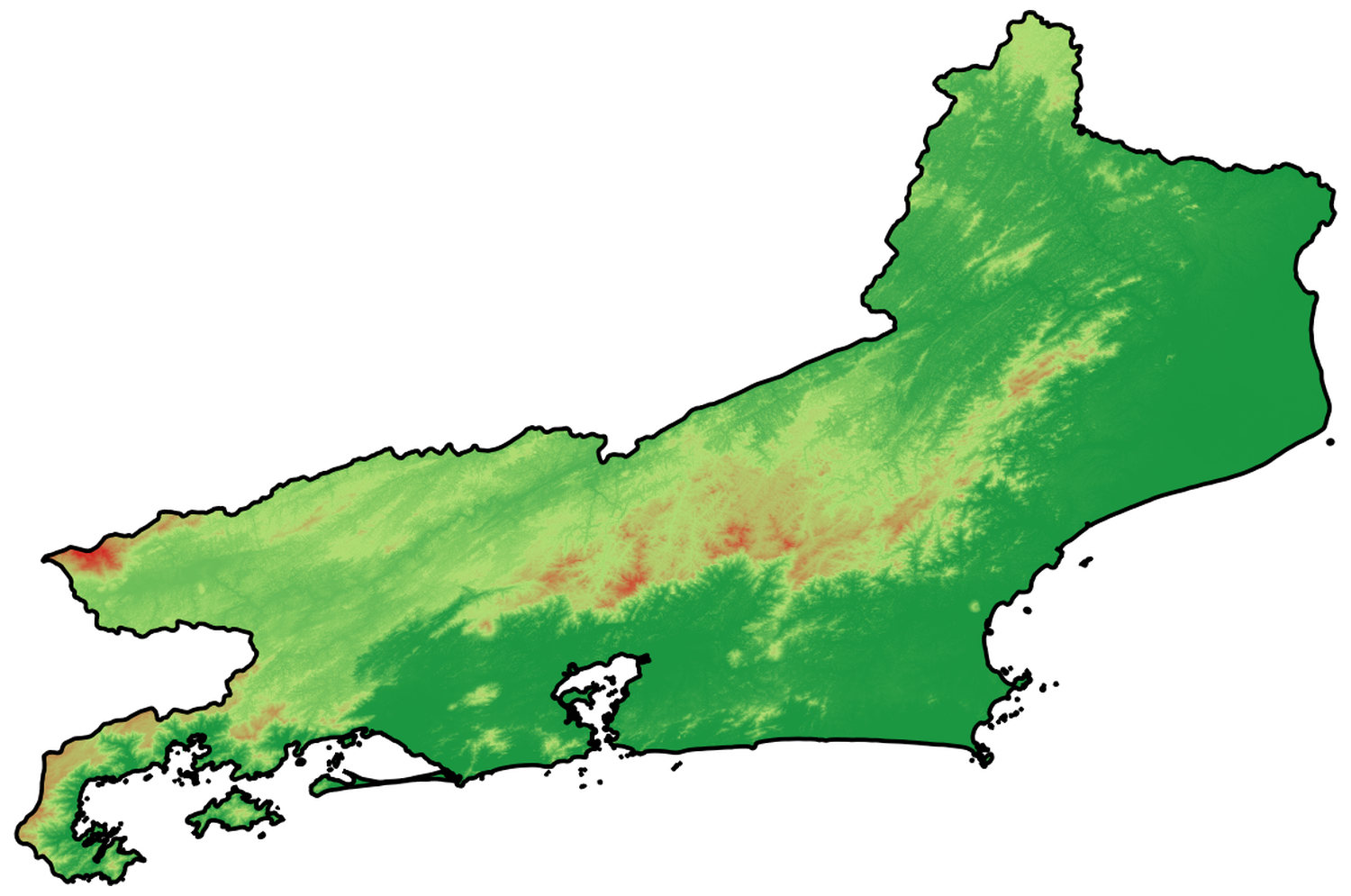}
    \caption*{\footnotesize Elevation Model (DEM)}\end{subfigure}
  \caption{Geospatial input datasets for Rio de Janeiro State evaluated in HERA-S.}
  \label{fig:fig4}
\end{figure*}

The methodology was applied across the State of Rio de Janeiro, Brazil. Input data included state boundaries, transport routes, transmission lines, substations, SRTM 90\,m DEM elevation data, regional geology layers, and federal/state conservation areas (Figure~\ref{fig:fig4}).

\subsection{Sequential Filtering}
Let $\Omega$ represent all DEM grid cells covering the Para\'iba do Sul river basin within Rio de Janeiro State. 

First, a topographic filter identified closed-loop pairs with gross head $>400$\,m and length-to-head ratio ($L:H$) $\le 10:1$, narrowing the search space to subset $\Omega_1 \subset \Omega$.

Subsequent filters removed cells excessively distant from high-voltage substations or major transport corridors to limit connection costs ($\Omega_2 \subset \Omega_1$). Protected environmental reserves were screened out ($\Omega_3 \subset \Omega_2$), followed by unfavorable geological formations ($\Omega_4 \subset \Omega_3$). Finally, sites with insufficient local watershed area were eliminated to ensure adequate initial water filling capability, producing final subset $\Omega_5 \subset \Omega_4$.

\subsection{Intensive Search Results}
\begin{figure*}[htbp]
  \centering
  \includegraphics[width=0.9\textwidth]{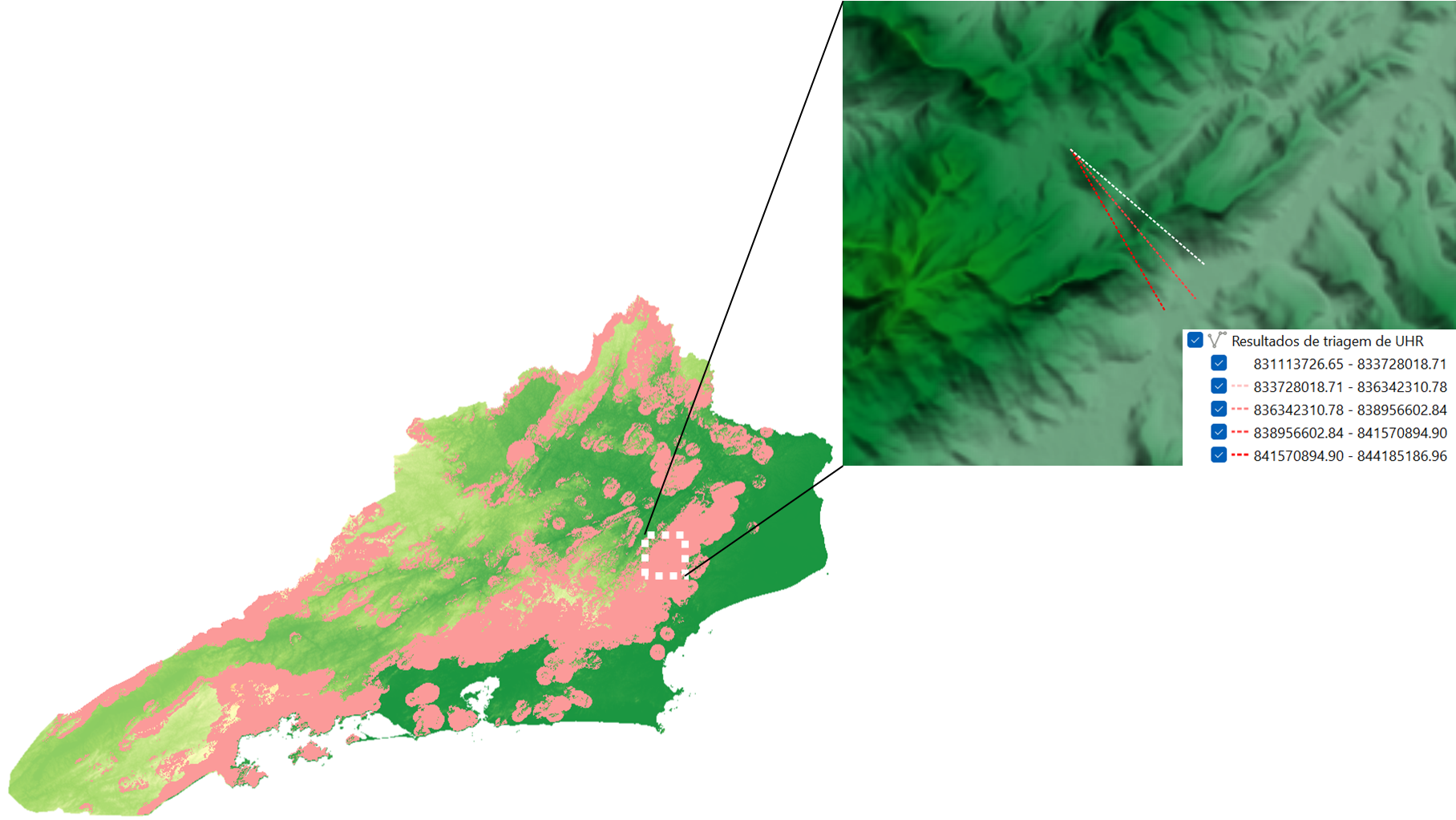}
  \caption{Intensive search candidate footprints (1000\,MW / 12\,h) in subset $\Omega_5$.}
  \label{fig:fig6}
\end{figure*}

Using subset $\Omega_5$, the intensive search evaluated 1000\,MW closed-loop configurations across storage durations from 6\,h to 120\,h. Figure~\ref{fig:fig6} shows candidate sites identified for 12\,h duration systems (dashed outlines) around a selected high-potential zone.

\subsection{Layout Optimization \& Cost Savings}
The reservoir optimization model was executed for a representative 1000\,MW / 12\,h (12,000\,MWh) project requiring $\approx 14.4\text{ million m}^3$ of storage volume under a gross head of 424\,m. The resulting mixed-integer program comprised over 13,000 decision variables (8,758 binary) and nearly 48,000 constraints.

\begin{figure*}[htbp]
  \centering
  \includegraphics[width=0.98\textwidth]{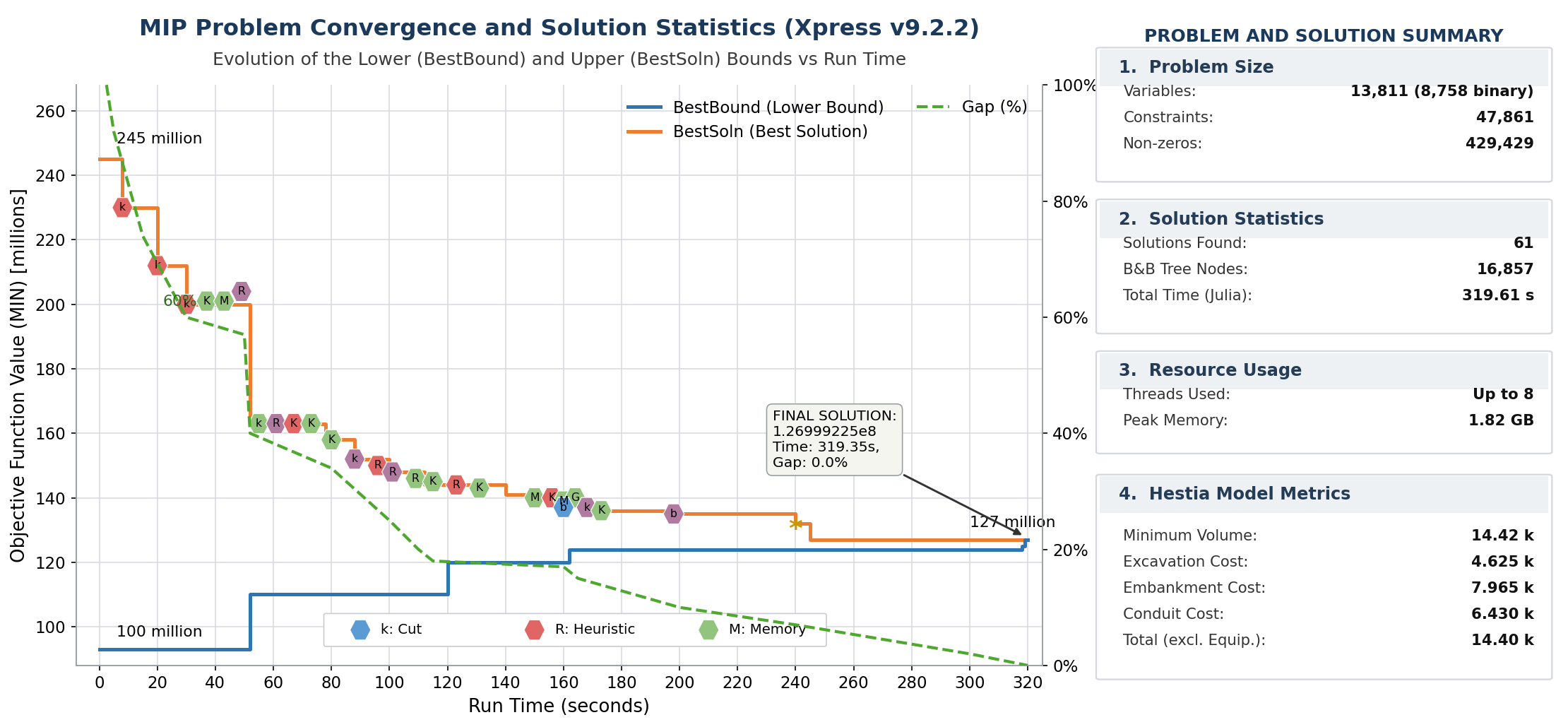}
  \caption{Optimization convergence curve showing upper bound (BestSoln), lower bound (BestBound), and optimality gap over time.}
  \label{fig:fig7}
\end{figure*}

Using the FICO Xpress solver, the model applied Branch-and-Bound with cut generation and aggressive heuristics. Figure~\ref{fig:fig7} shows convergence: an initial optimality gap exceeding 60\% was reduced to near zero in 320 seconds, converging to an optimized reservoir cost of USD 127 million.

The optimization prioritized balanced earthworks (embankment structures) over deep rock excavation. The solver discovered 61 feasible integer solutions during execution, reflecting the efficacy of the modified Traveling Salesman connectivity constraints.

Table~\ref{tab:costs} compares preliminary costs derived from geomorphon heuristics against the MIP-optimized design. The optimization reduced total project CAPEX by USD 47 million ($\approx 6\%$).

\begin{table}[htbp]
  \centering
  \caption{Capital cost breakdown: Initial heuristic layout vs. MIP-optimized layout (USD million).}
  \label{tab:costs}
  \small
  \begin{tabular}{@{}lcc@{}}
    \toprule
    \textbf{Cost Element} & \textbf{Initial Heuristic} & \textbf{MIP Optimized}\\
    \midrule
    Upper reservoir            & 89  & 78 \\
    Lower reservoir            & 25  & 8  \\
    Tunnels and canals         & 91  & 91 \\
    Powerhouse                 & 14  & 14 \\
    Electromechanical equipment& 226 & 226\\
    Grid connection            & 20  & 20 \\
    Access roads               & 4   & 4  \\
    Socio-environmental        & 47  & 44 \\
    Indirects \& contingencies & 258 & 242\\
    \midrule
    \textbf{Total Cost}         & \textbf{774} & \textbf{727}\\
    \textbf{Cost Reduction}    & \multicolumn{2}{c}{\textbf{6.1\%}}\\
    \bottomrule
  \end{tabular}
\end{table}

\begin{figure}[htbp]
  \centering
  \includegraphics[width=0.78\linewidth]{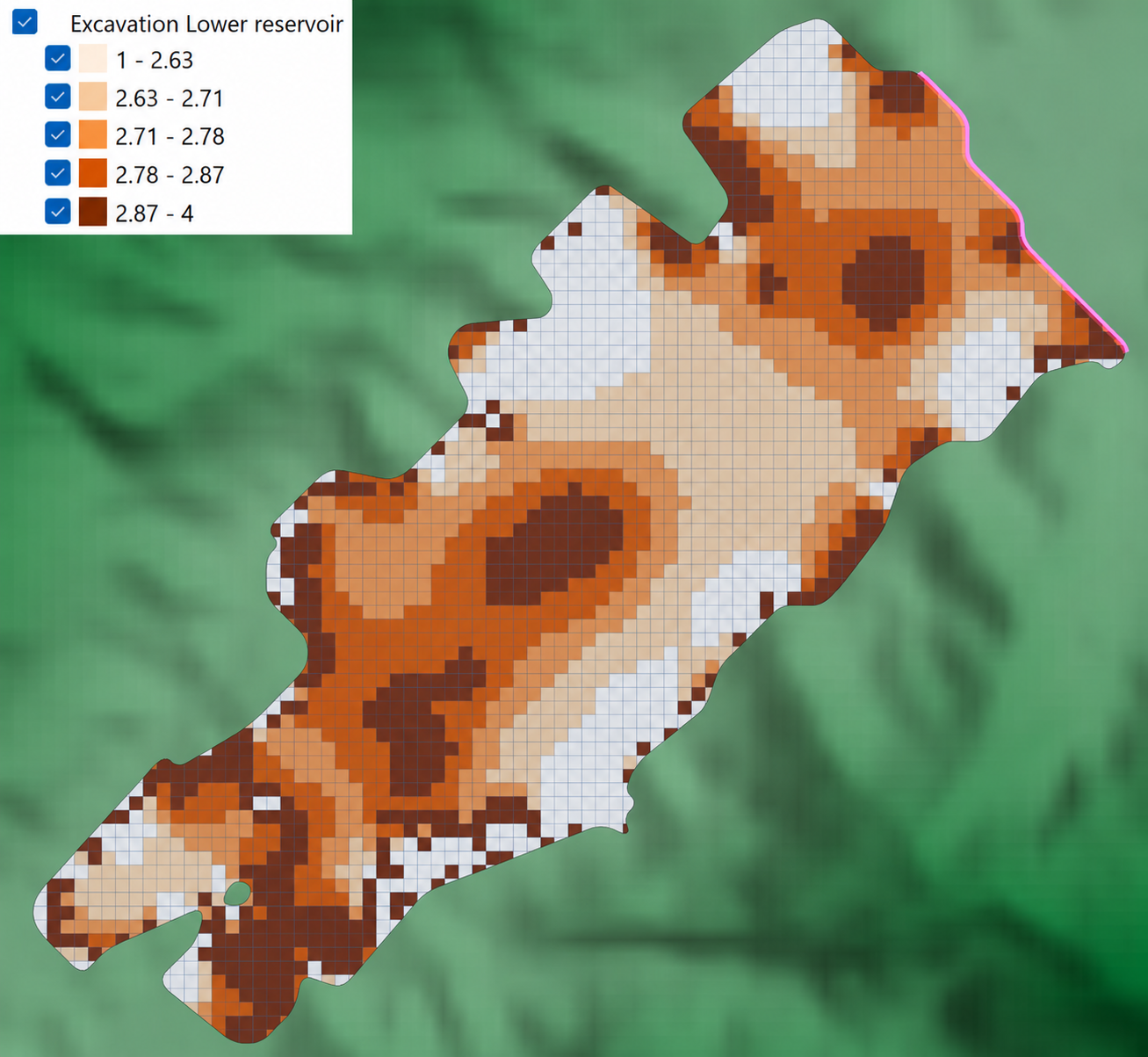}
  \caption{Optimized lower reservoir earthwork excavation depth map.}
  \label{fig:fig9}
\end{figure}

Significant savings occurred in the lower reservoir works, where the MIP model utilized local excavation material (up to 4\,m depth) to construct the retaining embankments directly (Figure~\ref{fig:fig9}).

\subsection{Engineering Layout}
\begin{figure}[htbp]
  \centering
  \includegraphics[width=0.82\linewidth]{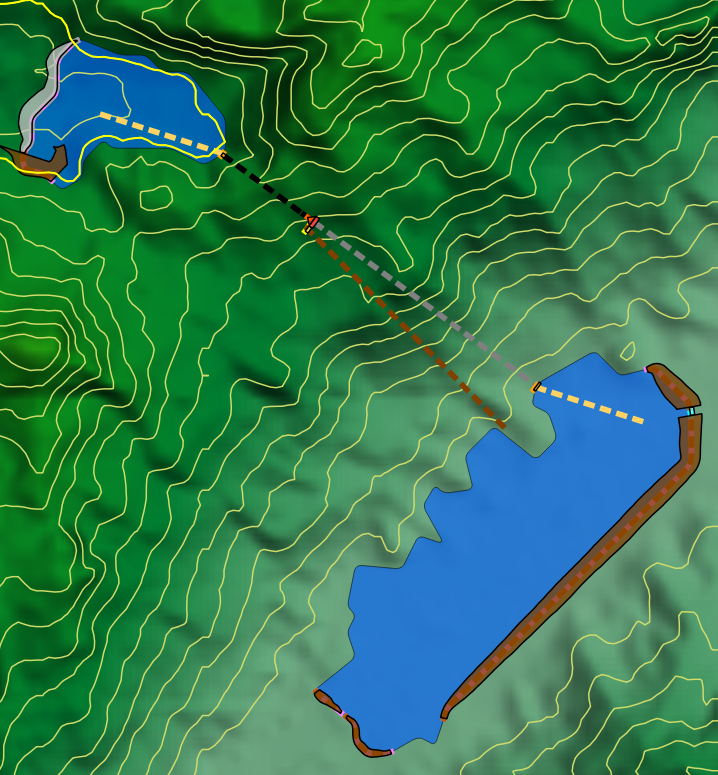}
  \caption{3D engineering layout generated by HERA-S for the 1000\,MW / 12\,h PSH plant in S\~ao Fid\'elis.}
  \label{fig:fig11}
\end{figure}

For top-ranking sites, HERA-S generated detailed engineering drawings and itemized bills of quantities. Figure~\ref{fig:fig11} shows the primary structures for a selected site in S\~ao Fid\'elis municipality: concrete/earth dams (grey/brown), spillway (light blue), reservoirs (dark blue), intake structures (orange), powerhouse (red), access tunnel (dashed brown), water canals (dashed yellow), and penstock tunnel (grey/black dashed lines).

The project achieves a unit storage cost of $\approx 64$\,USD/kWh ($727$\,million USD total CAPEX), significantly lower than equivalent BESS utility-scale deployments for 12\,h discharge durations.

\section{Conclusions}
\label{sec:concl}

This paper presented HERA-S, an integrated decision-support tool for prospecting and optimizing closed-loop pumped-storage hydropower (PSH) projects. By coupling GIS spatial filtering with mixed-integer mathematical programming, HERA-S automates site screening, reservoir shape optimization, and civil engineering cost estimation.

The core optimization formulation uses binary decision variables per DEM cell to balance excavation and fill earthworks while enforcing topological connectivity via separating planes and subtour elimination constraints. Applied to a 1000\,MW / 12,000\,MWh case study in Rio de Janeiro State, the optimization converged within 320 seconds, cutting total project CAPEX by 6\% (USD 47 million) compared to standard GIS heuristics.

The resulting unit storage cost ($\approx 64$\,USD/kWh) demonstrates the clear cost advantage of long-duration PSH over lithium-ion batteries.

Future extensions of HERA-S include: (i) full integration with expansion tools like OptGen; (ii) life-cycle carbon assessments comparing PSH with chemical BESS; (iii) explicit modeling of transmission congestion relief benefits; and (iv) expanding site inventories across additional mountainous regions in Brazil.


\section*{About the Authors}
\small
\begin{description}[leftmargin=0pt,itemsep=0.5em]
\item[\textbf{Luiz Rodolpho ALBUQUERQUE}] Expert at PSR since 2015; M.Sc. in Urban and Environmental Engineering (PUC-Rio), B.Sc. in Civil Engineering (UFRJ). Lead developer of the HERA engineering module for conventional and pumped-storage hydropower.

\item[\textbf{Rafael KELMAN}] Executive Director at PSR since 1997; Ph.D. in Systems Engineering (COPPE/UFRJ), M.Sc. in Water Resources. International consultant to the World Bank and IDB across 30+ countries.

\item[\textbf{Tarc\'isio CASTRO}] Team Leader at PSR since 2005; M.Sc. in Urban and Environmental Engineering (PUC-Rio). Assistant Professor at Poli/UFRJ since 1980, specializing in hydrological modeling and environmental due diligence.

\item[\textbf{Ana PETRUNGARO}] Lead Analyst at PSR since 2025; M.Sc. candidate in Prof\'Agua (UERJ), MBA (Univ. of Akron). Specialist in dam safety auditing, flood risk mapping, and water resources planning.

\item[\textbf{Tiago ANDRADE}] Lead Specialist at PSR since 2018; Ph.D. in Production Engineering (PUC-Rio/RMIT). Optimization algorithm developer focusing on linear, stochastic, and decomposition methods for SDDP and HERA.
\end{description}

\balance
\end{document}